# Algorithm to compute the rank and a Cartan subalgebra of a matrix Lie algebra with *Mathematica*

**Pablo Alberca Bjerregaard**

*Department of Applied Mathematics*
*Escuela Técnica Superior de Ingeniería Industrial*
*Campus de El Ejido*
*Universidad de Málaga*
*Málaga –Spain*
*E–Mail:pgalberca@uma.es*

**Cándido Martín González**

*Department of Algebra, Geometry and Topology*
*Facultad de Ciencias*
*Universidad de Málaga*
*Campus de Teatinos*
*Málaga –Spain*
*E–Mail:candido@apncs.cie.uma.es*

*2001.*

**Abstract.** *We present in this paper a set of routines constructed to compute the rank of a matrix Lie algebra and also to determine a Cartan subalgebra from a given list of elements.*

# ■ Introduction

It is clear the usefullness of the Cartan subalgebras in the Lie algebra context (see the references). We have worked in a previous paper in the links between the Lie algebras and particle physics. In that work we needed a routine to construct Cartan subalgebras (in order to define observables) and, as a consequence, the rank of a Lie algebra. In this paper we present all the routines defined for that purpose.

We will show some examples with clasical Lie algebras (like *so(n)*, *sl(n)*) and with two models in particle physics (*split* $\mathfrak{g}_2$ y $so(4, 4) \oplus so(2, 2)$). The theoretic result on which the rountine is based is the following:

## ■ Teorema.

*Let* $\mathfrak{L}$ *be a real semisimple matrix Lie algebra finite−dimensional.Let* $\mathfrak{H}$ *be a subalgebra of* $\mathfrak{L}$. *Then* $\mathfrak{H}$ *is a Cartan subalgebra if, and only if, it is abelian maximal and for all x and y in* $\mathfrak{H}$, $[x, y^t] = 0$.

The routine proceeds as follows: first of all a list of elements of the algebra is introduced. Once the necessary conditions are verified, its maximality is analyzed. If it is maximal, they form a Cartan subalgebra and its cardinal is the rank of the Lie algebra, but if they are not, the routine solve the system for the possible elements that can extend the set in order to form a Cartan subalgebra. This process is repeated until the subalgebra is maximal.

# ■ Auxiliary commands

The following command computes the intersection of two subspaces. :

```
δi_,j_ := If[i == j, 1, 0];
Cap[m_, n_] := Module[{},
gen1 := Sum[x_i m[[i]], {i, 1, Length[m]}]; gen2 := Sum[y_i n[[i]], {i, 1, Length[n]}];
in = gen1 //. ToRules[Reduce[gen1 == gen2]];
varia := Variables[in]; num := Length[varia];
rules := Table[varia[[j]] → δi,j , {i, 1, num}, {j, 1, num}];
Table[in //. rules[[i]], {i, 1, num}]];
```

and of course the Lie bracket:

```
c[x_, y_] := x.y - y.x;
```

We start from a matrix Lie algebra $\mathfrak{L}$ with basis $\mathfrak{B}$ constituted by the matrices $\{b_j\}_{j=1}^n$ (previously defined). The future basis of the Cartan subalgebra $\mathcal{H}$ will be denoted by $\{v_j\}_{j=1}^r$ where $r$ is the rank. The next function, **KerAd[m,b,n]** where **m** is a matrix, **b** is the notation for the basis and **n** is the dimension of the algebra, computes de kernel of the linear map:

$$\begin{array}{rccc} \mathrm{ad}_{\mathrm{m}}: & \mathcal{L} & \longrightarrow & \mathcal{L} \\ & \mathrm{x} & \longmapsto & [\mathrm{m}, \mathrm{x}] \end{array} \tag{1}$$

where [m, x] is the Lie bracket in $\mathcal{L}$, that is to say, [m,x]=mx−xm:

```
KerAd[m_, b_, n_] := Module[{},
X = Sum[a_i b_i, {i, 1, n}]; alla = Table[a_t , {t, 1, n}];
Matr = {};
Do[AppendTo[Matr, alla //. Solve[X == c[m, b_l]][[1]]], {l,
Matr = Transpose[Matr]; NullSpace[Matr]]
```

The next module defines **Indepen[v,r]**, which checks if the elements are linearly independent:

```
Indepen[v_, r_] := Module[{fal = 1, indep = {}},
indep = ToRules[Reduce[Sum[d_k v_k, {k, 1, r}] == θ]];
  Do[If[((d_k //. indep) === 0) == False, fal = 0], {k, 1, 
```

The next one, **Co[v,r]**, check if the set is abelian:

```
Co[v_, r_] := Module[{fal = 1},
Do[If[c[vi , vj] ≠ θ, fal = 0], {i, 1, r}, {j, i + 1, r}]; fa
```

The module **Star[v,r]** check if the set $\{v_j\}_{j=1}^{r}$ satisfies the condition $[x, x^t] = 0$ for all $x$:

```
Star[v_, r_] := Module[{fal = 1},
Do[If[c[vi , Transpose[vj]] ≠ θ, fal = 0],
            {i, 1, r}, {j, i, r}]; fal]
```

Now, the command is **Maximal[v,r,b,n]** which establishes when an abelian set is maximal. :

```
Maximal[v_, r_, b_, n_] := Module[{fal = 1},
ec = {c[X, Transpose[X]] == θ};
Do[
   AppendTo[ec, c[vk , X] == θ], {k, 1, r}];
        Do[AppendTo[ec,
    c[Transpose[vk], X] == θ], {k, 1, r}];
                Do[AppendTo[ec,
    Scala[X, vk] == 0], {k, 1, r}];
rulesol = Reduce[ec];
  If[Length[rulesol] ≠ n, fal = 0];
        fal]
```

At last, we have the module **Addone[r,v,b,n]** which adds, inside the construction process of a Cartan subalgebra, a new element in each step:

```
Addone[r_, v_, b_, n_] := Module[{},
verif := c[Y, Transpose[Y]] == θ;
Do[verif = verif && c[Y, vj] == θ &&
                    c[Y, Transpose[vj]] == θ && Scala[Y
    {j, r - 1}];
ta = Y //. {ToRules[Reduce[verif]]}[[1]];
        var = Variables[ta];
num = Length[var];
        rul = Table[var[[j]] → δi,j , {i, num}, {j, num}];
        kk = 1;
While[(ta //. rul[[kk]]) == θ , kk = kk + 1];
        vr = ta //. rul[[kk]]]
```

## ■ Main routine

It is time to present the main routine. We need first to introduce a basis of the particular Lie algebra we are interested in. We are going to present a comand for each algebra but it is important to say that, of course, the routine remains the same.

## ■ Lie algebra so(*n*)

We compute first a basis of the Lie algebra:

```
BSo[t_] := Module[{kk = 1},
n = 1/2 (t - 1) t; θ = Table[0, {i, 1, t}, {j, 1, t}];
Do[m = θ; m[[p, q]] = 1; F_{p,q} = m, {q, 1, t}, {p, 1, t}];
Do[J_{p,q} = F_{p,q} - F_{q,p} , {p, 1, t}, {q, 1, t}]; kk = 1;
Do[Do[b_kk = J_{p,q}; kk = kk + 1, {q, p + 1, t}], {p, 1, t - 1}];
Scala[x_, y_] := Sum_{i=1}^{t} Sum_{j=1}^{t} x[[i, j]] y[[i, j]]; X = Sum_{i=1}^{n} a_i b_i ;]
```

and the the main routine is:

```
RankOfSo[t_, dat_] := Module[{},
r = Length[dat];
Do[v_k = dat[[k]], {k, 1, r}];
    If[Indepen[v, r] == 1 && Co[v, r] == 1 && Star[v, r] ==
        pset = KerAd[v_1 , b, n]; If[r > 1,
                        Do[pset = Cap[pset, KerAd[v_i , b, n]

Y = Sum_{k=1}^{Length[pset]} f_k Sum_{i=1}^{n} pset[[k, i]] b_i ;
While[Maximal[v, r, b, n] == 0, r = r + 1; Addone[r, v, b, ı
Print["The rank is ", r, ",", " and a Cartan subalgebr
Do[comb = Solve[v_q == X];
    Print["v"_q , " = ", Sum_{p=1}^{n} (a_p //. comb[[1]]) * b_p], {q, 1, ı
                    Print[False]]]
```

The fundamental line is

```
While[Maximal[v, r, b, n] == 0,
 r = r + 1; Addone[r, v, b, n]]
```

Let’s look an example:

```
BSo[6]
```

and then a generic element is

**X**

$$\begin{pmatrix} 0 & a_1 & a_2 & a_3 & a_4 & a_5 \\ -a_1 & 0 & a_6 & a_7 & a_8 & a_9 \\ -a_2 & -a_6 & 0 & a_{10} & a_{11} & a_{12} \\ -a_3 & -a_7 & -a_{10} & 0 & a_{13} & a_{14} \\ -a_4 & -a_8 & -a_{11} & -a_{13} & 0 & a_{15} \\ -a_5 & -a_9 & -a_{12} & -a_{14} & -a_{15} & 0 \end{pmatrix}$$

We make then

```
RankOfSo[6, {b1}]
```

The rank is 3, and a Cartan subalgebra:

$v_1 = b_1$

$v_2 = b_{15}$

$v_3 = b_{10}$

where we have also shown the time it took. We have then that the rank is 3 and a Cartan subalgebra is constituted by the matrices

```
{{b1, b10, b15}}
```

$$\left( \begin{pmatrix} 0 & 1 & 0 & 0 & 0 & 0 \\ -1 & 0 & 0 & 0 & 0 & 0 \\ 0 & 0 & 0 & 0 & 0 & 0 \\ 0 & 0 & 0 & 0 & 0 & 0 \\ 0 & 0 & 0 & 0 & 0 & 0 \\ 0 & 0 & 0 & 0 & 0 & 0 \end{pmatrix} \begin{pmatrix} 0 & 0 & 0 & 0 & 0 & 0 \\ 0 & 0 & 0 & 0 & 0 & 0 \\ 0 & 0 & 0 & 1 & 0 & 0 \\ 0 & 0 & -1 & 0 & 0 & 0 \\ 0 & 0 & 0 & 0 & 0 & 0 \\ 0 & 0 & 0 & 0 & 0 & 0 \end{pmatrix} \begin{pmatrix} 0 & 0 & 0 & 0 & 0 & 0 \\ 0 & 0 & 0 & 0 & 0 & 0 \\ 0 & 0 & 0 & 0 & 0 & 0 \\ 0 & 0 & 0 & 0 & 0 & 0 \\ 0 & 0 & 0 & 0 & 0 & 1 \\ 0 & 0 & 0 & 0 & -1 & 0 \end{pmatrix} \right)$$

We can start form two elements, $b_1$ y $b_{10}$ , and we have

```
RankOfSo[6, {b1, b10}]
```

The rank is 3, and a Cartan subalgebra:

$v_1 = b_1$

$v_2 = b_{10}$

$v_3 = b_{15}$

Another example could be

```
BSo[9]
```

and the a generic element is

```
X
```

$$\begin{pmatrix} 0 & a_1 & a_2 & a_3 & a_4 & a_5 & a_6 & a_7 & a_8 \\ -a_1 & 0 & a_9 & a_{10} & a_{11} & a_{12} & a_{13} & a_{14} & a_{15} \\ -a_2 & -a_9 & 0 & a_{16} & a_{17} & a_{18} & a_{19} & a_{20} & a_{21} \\ -a_3 & -a_{10} & -a_{16} & 0 & a_{22} & a_{23} & a_{24} & a_{25} & a_{26} \\ -a_4 & -a_{11} & -a_{17} & -a_{22} & 0 & a_{27} & a_{28} & a_{29} & a_{30} \\ -a_5 & -a_{12} & -a_{18} & -a_{23} & -a_{27} & 0 & a_{31} & a_{32} & a_{33} \\ -a_6 & -a_{13} & -a_{19} & -a_{24} & -a_{28} & -a_{31} & 0 & a_{34} & a_{35} \\ -a_7 & -a_{14} & -a_{20} & -a_{25} & -a_{29} & -a_{32} & -a_{34} & 0 & a_{36} \\ -a_8 & -a_{15} & -a_{21} & -a_{26} & -a_{30} & -a_{33} & -a_{35} & -a_{36} & 0 \end{pmatrix}$$

We can make then

```
RankOfSo[9, {b1}]
```

The rank is 4, and a Cartan subalgebra:

$v_1 = b_1$

$v_2 = b_{36}$

$v_3 = b_{31}$

$v_4 = b_{22}$

As a final example, we can use the Lie algebra *so(32)*, recently used in string theory (dimension= $\frac{32}{2}$ 31 = 496):

```
BSo[32]
```

and a generic element

```
X
```

$$
\begin{pmatrix}
0 & a_{1} & a_{2} & a_{3} & a_{4} & a_{5} & a_{6} & a_{7} & a_{8} & a_{9} & a_{10} & a_{11} & a_{12} & a_{13} & a_{14} & a_{15} & a_{16} & a_{17} & a_{18} & a \\
-a_{1} & 0 & a_{32} & a_{33} & a_{34} & a_{35} & a_{36} & a_{37} & a_{38} & a_{39} & a_{40} & a_{41} & a_{42} & a_{43} & a_{44} & a_{45} & a_{46} & a_{47} & a_{48} & a \\
-a_{2} & -a_{32} & 0 & a_{62} & a_{63} & a_{64} & a_{65} & a_{66} & a_{67} & a_{68} & a_{69} & a_{70} & a_{71} & a_{72} & a_{73} & a_{74} & a_{75} & a_{76} & a_{77} & a \\
-a_{3} & -a_{33} & -a_{62} & 0 & a_{91} & a_{92} & a_{93} & a_{94} & a_{95} & a_{96} & a_{97} & a_{98} & a_{99} & a_{100} & a_{101} & a_{102} & a_{103} & a_{104} & a_{105} & a \\
-a_{4} & -a_{34} & -a_{63} & -a_{91} & 0 & a_{119} & a_{120} & a_{121} & a_{122} & a_{123} & a_{124} & a_{125} & a_{126} & a_{127} & a_{128} & a_{129} & a_{130} & a_{131} & a_{132} & a \\
-a_{5} & -a_{35} & -a_{64} & -a_{92} & -a_{119} & 0 & a_{146} & a_{147} & a_{148} & a_{149} & a_{150} & a_{151} & a_{152} & a_{153} & a_{154} & a_{155} & a_{156} & a_{157} & a_{158} & a \\
-a_{6} & -a_{36} & -a_{65} & -a_{93} & -a_{120} & -a_{146} & 0 & a_{172} & a_{173} & a_{174} & a_{175} & a_{176} & a_{177} & a_{178} & a_{179} & a_{180} & a_{181} & a_{182} & a_{183} & a \\
-a_{7} & -a_{37} & -a_{66} & -a_{94} & -a_{121} & -a_{147} & -a_{172} & 0 & a_{197} & a_{198} & a_{199} & a_{200} & a_{201} & a_{202} & a_{203} & a_{204} & a_{205} & a_{206} & a_{207} & a \\
-a_{8} & -a_{38} & -a_{67} & -a_{95} & -a_{122} & -a_{148} & -a_{173} & -a_{197} & 0 & a_{221} & a_{222} & a_{223} & a_{224} & a_{225} & a_{226} & a_{227} & a_{228} & a_{229} & a_{230} & a \\
-a_{9} & -a_{39} & -a_{68} & -a_{96} & -a_{123} & -a_{149} & -a_{174} & -a_{198} & -a_{221} & 0 & a_{244} & a_{245} & a_{246} & a_{247} & a_{248} & a_{249} & a_{250} & a_{251} & a_{252} & a \\
-a_{10} & -a_{40} & -a_{69} & -a_{97} & -a_{124} & -a_{150} & -a_{175} & -a_{199} & -a_{222} & -a_{244} & 0 & a_{266} & a_{267} & a_{268} & a_{269} & a_{270} & a_{271} & a_{272} & a_{273} & a \\
-a_{11} & -a_{41} & -a_{70} & -a_{98} & -a_{125} & -a_{151} & -a_{176} & -a_{200} & -a_{223} & -a_{245} & -a_{266} & 0 & a_{287} & a_{288} & a_{289} & a_{290} & a_{291} & a_{292} & a_{293} & a \\
-a_{12} & -a_{42} & -a_{71} & -a_{99} & -a_{126} & -a_{152} & -a_{177} & -a_{201} & -a_{224} & -a_{246} & -a_{267} & -a_{287} & 0 & a_{307} & a_{308} & a_{309} & a_{310} & a_{311} & a_{312} & a \\
-a_{13} & -a_{43} & -a_{72} & -a_{100} & -a_{127} & -a_{153} & -a_{178} & -a_{202} & -a_{225} & -a_{247} & -a_{268} & -a_{288} & -a_{307} & 0 & a_{326} & a_{327} & a_{328} & a_{329} & a_{330} & a \\
-a_{14} & -a_{44} & -a_{73} & -a_{101} & -a_{128} & -a_{154} & -a_{179} & -a_{203} & -a_{226} & -a_{248} & -a_{269} & -a_{289} & -a_{308} & -a_{326} & 0 & a_{344} & a_{345} & a_{346} & a_{347} & a \\
-a_{15} & -a_{45} & -a_{74} & -a_{102} & -a_{129} & -a_{155} & -a_{180} & -a_{204} & -a_{227} & -a_{249} & -a_{270} & -a_{290} & -a_{309} & -a_{327} & -a_{344} & 0 & a_{361} & a_{362} & a_{363} & a \\
-a_{16} & -a_{46} & -a_{75} & -a_{103} & -a_{130} & -a_{156} & -a_{181} & -a_{205} & -a_{228} & -a_{250} & -a_{271} & -a_{291} & -a_{310} & -a_{328} & -a_{345} & -a_{361} & 0 & a_{377} & a_{378} & a \\
-a_{17} & -a_{47} & -a_{76} & -a_{104} & -a_{131} & -a_{157} & -a_{182} & -a_{206} & -a_{229} & -a_{251} & -a_{272} & -a_{292} & -a_{311} & -a_{329} & -a_{346} & -a_{362} & -a_{377} & 0 & a_{392} & a \\
-a_{18} & -a_{48} & -a_{77} & -a_{105} & -a_{132} & -a_{158} & -a_{183} & -a_{207} & -a_{230} & -a_{252} & -a_{273} & -a_{293} & -a_{312} & -a_{330} & -a_{347} & -a_{363} & -a_{378} & -a_{392} & 0 & a \\
-a_{19} & -a_{49} & -a_{78} & -a_{106} & -a_{133} & -a_{159} & -a_{184} & -a_{208} & -a_{231} & -a_{253} & -a_{274} & -a_{294} & -a_{313} & -a_{331} & -a_{348} & -a_{364} & -a_{379} & -a_{393} & -a_{406} & \\
-a_{20} & -a_{50} & -a_{79} & -a_{107} & -a_{134} & -a_{160} & -a_{185} & -a_{209} & -a_{232} & -a_{254} & -a_{275} & -a_{295} & -a_{314} & -a_{332} & -a_{349} & -a_{365} & -a_{380} & -a_{394} & -a_{407} & -a \\
-a_{21} & -a_{51} & -a_{80} & -a_{108} & -a_{135} & -a_{161} & -a_{186} & -a_{210} & -a_{233} & -a_{255} & -a_{276} & -a_{296} & -a_{315} & -a_{333} & -a_{350} & -a_{366} & -a_{381} & -a_{395} & -a_{408} & -a \\
-a_{22} & -a_{52} & -a_{81} & -a_{109} & -a_{136} & -a_{162} & -a_{187} & -a_{211} & -a_{234} & -a_{256} & -a_{277} & -a_{297} & -a_{316} & -a_{334} & -a_{351} & -a_{367} & -a_{382} & -a_{396} & -a_{409} & -a \\
-a_{23} & -a_{53} & -a_{82} & -a_{110} & -a_{137} & -a_{163} & -a_{188} & -a_{212} & -a_{235} & -a_{257} & -a_{278} & -a_{298} & -a_{317} & -a_{335} & -a_{352} & -a_{368} & -a_{383} & -a_{397} & -a_{410} & -a \\
-a_{24} & -a_{54} & -a_{83} & -a_{111} & -a_{138} & -a_{164} & -a_{189} & -a_{213} & -a_{236} & -a_{258} & -a_{279} & -a_{299} & -a_{318} & -a_{336} & -a_{353} & -a_{369} & -a_{384} & -a_{398} & -a_{411} & -a \\
-a_{25} & -a_{55} & -a_{84} & -a_{112} & -a_{139} & -a_{165} & -a_{190} & -a_{214} & -a_{237} & -a_{259} & -a_{280} & -a_{300} & -a_{319} & -a_{337} & -a_{354} & -a_{370} & -a_{385} & -a_{399} & -a_{412} & -a \\
-a_{26} & -a_{56} & -a_{85} & -a_{113} & -a_{140} & -a_{166} & -a_{191} & -a_{215} & -a_{238} & -a_{260} & -a_{281} & -a_{301} & -a_{320} & -a_{338} & -a_{355} & -a_{371} & -a_{386} & -a_{400} & -a_{413} & -a \\
-a_{27} & -a_{57} & -a_{86} & -a_{114} & -a_{141} & -a_{167} & -a_{192} & -a_{216} & -a_{239} & -a_{261} & -a_{282} & -a_{302} & -a_{321} & -a_{339} & -a_{356} & -a_{372} & -a_{387} & -a_{401} & -a_{414} & -a \\
-a_{28} & -a_{58} & -a_{87} & -a_{115} & -a_{142} & -a_{168} & -a_{193} & -a_{217} & -a_{240} & -a_{262} & -a_{283} & -a_{303} & -a_{322} & -a_{340} & -a_{357} & -a_{373} & -a_{388} & -a_{402} & -a_{415} & -a \\
-a_{29} & -a_{59} & -a_{88} & -a_{116} & -a_{143} & -a_{169} & -a_{194} & -a_{218} & -a_{241} & -a_{263} & -a_{284} & -a_{304} & -a_{323} & -a_{341} & -a_{358} & -a_{374} & -a_{389} & -a_{403} & -a_{416} & -a \\
-a_{30} & -a_{60} & -a_{89} & -a_{117} & -a_{144} & -a_{170} & -a_{195} & -a_{219} & -a_{242} & -a_{264} & -a_{285} & -a_{305} & -a_{324} & -a_{342} & -a_{359} & -a_{375} & -a_{390} & -a_{404} & -a_{417} & -a \\
-a_{31} & -a_{61} & -a_{90} & -a_{118} & -a_{145} & -a_{171} & -a_{196} & -a_{220} & -a_{243} & -a_{265} & -a_{286} & -a_{306} & -a_{325} & -a_{343} & -a_{360} & -a_{376} & -a_{391} & -a_{405} & -a_{418} & -a
\end{pmatrix}
$$

```
RankOfSo[32, {b1 }]
```

The rank is 16, and a Cartan subalgebra:

$v_1 = b_1$

$v_2 = b_{496}$

$v_3 = b_{491}$

$v_4 = b_{482}$

$v_5 = b_{469}$

$v_6 = b_{452}$

$v_7 = b_{431}$

$v_8 = b_{406}$

$v_9 = b_{377}$

$v_{10} = b_{344}$

$v_{11} = b_{307}$

$v_{12} = b_{266}$

$v_{13} = b_{221}$

$v_{14} = b_{172}$

$v_{15} = b_{119}$

$v_{16} = b_{62}$

## ■ Lie algebra sl(*n*)

The second example could be the Lie algebra *sl(n)*. First of all we compute a basis as follows

```
Bsl[t_] := Module[{k = 1}, n = t^2 - 1;
θ = Table[0, {i, 1, t}, {j, 1, t}];
Do[m = θ; m[[p, q]] = 1; F_{p,q} = m, {q, 1, t}, {p, 1, t}];
Do[Do[If[p ≠ q, b_k = F_{p,q}; k = k + 1], {q, 1, t}], {p, 1, t}]
Do[b_k = F_{1,1} - F_{p,p}; k = k + 1, {p, 2, t}];
Scala[x_, y_] := Sum[Sum[x[[i, j]] y[[i, j]], {j, 1, t}], {i, 1, t}]; X = Sum[a_i b_i, {i, 1, n}]; ]
```

Then we have

```
RankOfsl[t_, dat_] := Module[{},
r = Length[dat];
Do[v_k = dat[[k]], {k, 1, r}];
    If[Indepen[v, r] == 1 && Co[v, r] == 1 && Star[v, r] ==
        pset = KerAd[v_1, b, n]; If[r > 1,
                        Do[pset = Cap[pset, KerAd[v_i, b, n]
Y = Sum[f_k Sum[pset[[k, i]] b_i, {i, 1, n}], {k, 1, Length[pset]}];
While[Maximal[v, r, b, n] == 0, r = r + 1; Addone[r, v, b, ı
Print["The rank is ", r, ",", " and a Cartan subalgebr
Do[comb = Solve[v_q == X];
    Print["v"_q, " = ", Sum[(a_p //. comb[[1]]) * b_p, {p, 1, n}], {q, 1, ı
                Print[False]]]
```

We can have a basis of sl(3) by making

```
Bsl[3];
```

and then a generic element is

```
X
```

$$\begin{pmatrix} a_7 + a_8 & a_1 & a_2 \\ a_3 & -a_7 & a_4 \\ a_5 & a_6 & -a_8 \end{pmatrix}$$

If we take the first element of the basis, $b_1$ , we obtain

```
RankOfsl[3, {b1}]
```

False

This is so becasuse

```
c[b1, Transpose[b1]]
```

$$\begin{pmatrix} 1 & 0 & 0 \\ 0 & -1 & 0 \\ 0 & 0 & 0 \end{pmatrix}$$

and it must verify $[x, x^t] = 0$. Nevertheless, if we take $b_8$ :

```
RankOfsl[3, {b7}]
```

The rank is 2, and a Cartan subalgebra:

$v_1 = b_7$

$v_2 = b_8$

we get that the rank is 2 and a Cartan subalgebra is

```
{{b7, b8}}
```

$$\left( \begin{pmatrix} 1 & 0 & 0 \\ 0 & -1 & 0 \\ 0 & 0 & 0 \end{pmatrix} \begin{pmatrix} 1 & 0 & 0 \\ 0 & 0 & 0 \\ 0 & 0 & -1 \end{pmatrix} \right)$$

In *sl(5)* we make first

```
Bsl[5];
```

with generic element

```
X
```

$$\begin{pmatrix} a_{21}+a_{22}+a_{23}+a_{24} & a_1 & a_2 & a_3 & a_4 \\ a_5 & -a_{21} & a_6 & a_7 & a_8 \\ a_9 & a_{10} & -a_{22} & a_{11} & a_{12} \\ a_{13} & a_{14} & a_{15} & -a_{23} & a_{16} \\ a_{17} & a_{18} & a_{19} & a_{20} & -a_{24} \end{pmatrix}$$

Then we can make

```
RankOfsl[5, {b21 , b22 }]
```

The rank is 4, and a Cartan subalgebra:

$v_1 = b_{21}$

$v_2 = b_{22}$

$v_3 = b_{23} - b_{24}$

$v_4 = -\frac{2 b_{21}}{3} - \frac{2 b_{22}}{3} + b_{23} + b_{24}$

to include these two elements in a Cartan subalgebra.

## ■ Lie algebra $\mathfrak{g}_2$

We present now the Lie algebra $\mathfrak{g}_2$ . A basis of this Lie algebra is

```
Bg2 := Module[{},
n = 14; t = 8; θ = Table[0, {i, 1, t}, {j, 1, t}];
Do[m = θ; m[[p, q]] = 1; F_{p,q} = m, {q, 1, t}, {p, 1, t}];
Do[J_{p,q} = F_{p,q} - F_{q,p} , {p, 1, t}, {q, 1, t}];
b_1 = J_{2,3} + J_{6,7} ; b_2 = J_{2,4} + J_{6,8} ; b_3 = J_{2,5} + J_{7,4} ;
b_4 = J_{2,6} + J_{8,4} ; b_5 = J_{2,7} + J_{4,5} ; b_6 = J_{2,8} + J_{4,6} ;
b_7 = J_{3,4} + J_{7,8} ; b_8 = J_{3,5} + J_{4,6} ; b_9 = J_{3,6} + J_{5,4} ;
b_10 = J_{3,7} + J_{8,4} ; b_11 = J_{3,8} + J_{4,7} ; b_12 = J_{5,6} + J_{7,8} ;
b_13 = J_{5,7} + J_{8,6} ; b_14 = J_{5,8} + J_{6,7} ;
Scala[x_, y_] := Sum_{i=1}^{t} Sum_{j=1}^{t} x[[i, j]] y[[i, j]]; X = Sum_{i=1}^{n} a_i b_i ; ]
```

and then a generic element is, after making,

```
Bg2
```

as

```
X
```

$$\begin{pmatrix} 0 & 0 & 0 & 0 & 0 & 0 & 0 & 0 \\ 0 & 0 & a_1 & a_2 & a_3 & a_4 & a_5 & a_6 \\ 0 & -a_1 & 0 & a_7 & a_8 & a_9 & a_{10} & a_{11} \\ 0 & -a_2 & -a_7 & 0 & a_5-a_9 & a_6+a_8 & a_{11}-a_3 & -a_4-a_{10} \\ 0 & -a_3 & -a_8 & a_9-a_5 & 0 & a_{12} & a_{13} & a_{14} \\ 0 & -a_4 & -a_9 & -a_6-a_8 & -a_{12} & 0 & a_1+a_{14} & a_2-a_{13} \\ 0 & -a_5 & -a_{10} & a_3-a_{11} & -a_{13} & -a_1-a_{14} & 0 & a_7+a_{12} \\ 0 & -a_6 & -a_{11} & a_4+a_{10} & -a_{14} & a_{13}-a_2 & -a_7-a_{12} & 0 \end{pmatrix}$$

We define then

```
RankOfg2[dat_] := Module[{},
r = Length[dat];
Do[v_k = dat[[k]], {k, 1, r}];
If[Indepen[v, r] == 1 && Co[v, r] == 1 && Star[v, r] == 1,
        pset = KerAd[v_1, b, n];
                If[r > 1, Do[pset = Cap[pset, KerAd[v_i, 

Y = Sum_{k=1}^{Length[pset]} f_k Sum_{i=1}^{n} pset[[k, i]] b_i;

While[Maximal[v, r, b, n] == 0, r = r + 1; Addone[r, v, b, 
Print["The rank is ", r, ",", " and a Cartan subalgebr
Do[comb = Solve[v_q == X];
    Print[(HoldForm[v])_q, " = ",
            Sum[(a_p //. comb[[1]]) * (HoldForm[b])_p, {p,
                Print[False]]]
```

By using the first element of the basis we have

```
RankOfg2[{b1}]
```

The rank is 2, and a Cartan subalgebra:

$v_1 = b_1$

$v_2 = b_{14}$

and the rank is 2. We can also try with

```
RankOfg2[{5 b5 + b1}]
```

The rank is 2, and a Cartan subalgebra:

$v_1 = b_1 + 5\, b_5$

$v_2 = b_1 + 5\, b_9$

## ■ Lie algebra *split* $\mathfrak{g}_2$

We now deal with a particle physic model: the Lie algebra *split* $\mathfrak{g}_2$ . A basis of this Lie algebra is

```
BSplitg2 := Module[{},
n = 14; t = 8; θ = Table[0, {i, 1, t}, {j, 1, t}];
Do[m = θ; m[[p, q]] = 1; F_{p,q} = m, {q, 1, t}, {p, 1, t}];
Do[J_{p,q} = F_{p,q} - F_{q,p}; Z_{p,q} = F_{p,q} + F_{q,p}, {p, 1, t}, {q, 1, t}]
b1 = J_{2,3} - J_{7,6}; b2 = -J_{4,2} + J_{6,8}; b3 = Z_{2,5} + Z_{4,7};
b4 = Z_{2,6} - Z_{4,8}; b5 = Z_{2,7} - Z_{4,5}; b6 = Z_{2,8} + Z_{4,6};
b7 = J_{3,4} + J_{7,8}; b8 = Z_{3,5} - Z_{4,6}; b9 = Z_{3,6} + Z_{4,5};
b10 = Z_{3,7} - Z_{4,8}; b11 = Z_{3,8} + Z_{4,7}; b12 = -J_{6,5} - J_{7,8};
b13 = J_{5,7} + J_{6,8}; b14 = J_{5,8} + J_{7,6};
Scala[x_, y_] := Sum_{i=1}^{t} Sum_{j=1}^{t} x[[i, j]] y[[i, j]]; X = Sum_{i=1}^{n} a_i b_i;]
```

and we can make

```
BSplitg2
```

in order to have that a generic element is as follows:

```
X
```

$$\begin{pmatrix} 0 & 0 & 0 & 0 & 0 & 0 & 0 & 0 \\ 0 & 0 & a_1 & a_2 & a_3 & a_4 & a_5 & a_6 \\ 0 & -a_1 & 0 & a_7 & a_8 & a_9 & a_{10} & a_{11} \\ 0 & -a_2 & -a_7 & 0 & a_9 - a_5 & a_6 - a_8 & a_3 + a_{11} & -a_4 - a_{10} \\ 0 & a_3 & a_8 & a_9 - a_5 & 0 & a_{12} & a_{13} & a_{14} \\ 0 & a_4 & a_9 & a_6 - a_8 & -a_{12} & 0 & a_1 - a_{14} & a_2 + a_{13} \\ 0 & a_5 & a_{10} & a_3 + a_{11} & -a_{13} & a_{14} - a_1 & 0 & a_7 - a_{12} \\ 0 & a_6 & a_{11} & -a_4 - a_{10} & -a_{14} & -a_2 - a_{13} & a_{12} - a_7 & 0 \end{pmatrix}$$

We define then

```
RankOfSplitg2[dat_] := Module[{},
r = Length[dat];
Do[v_k = dat[[k]], {k, 1, r}];
    If[Indepen[v, r] == 1 && Co[v, r] == 1 && Star[v, r] ==
        pset = KerAd[v_1 , b, n];
            If[r > 1, Do[pset = Cap[pset, KerAd[v_i , b, n]
                {i, 2, r}]];
```

$$\texttt{Y} = \sum_{k=1}^{\texttt{Length[pset]}} \texttt{f}_k \sum_{i=1}^{n} \texttt{pset}[\![k, i]\!]\, \texttt{b}_i ;$$

```
While[Maximal[v, r, b, n] == 0, r = r + 1; Addone[r, v, b, 
Print["The rank is ", r, ",", " and a Cartan subalgebr
Do[comb = Solve[v_q == X];
    Print[(HoldForm[v])_q, " = ",
    Sum[(a_p //. comb[[1]]) * (HoldForm[b])_p, {p, 1, n}]],
                    {q, 1, r}],     Print[False]]]
```

to check that the rank of this Lie algebra is also 2, as we have

```
RankOfSplitg2[{b_1 + b_2}]
```

The rank is 2, and a Cartan subalgebra:

$v_1 = b_1 + b_2$

$v_2 = 2\, b_4 - b_{10} + 3\, b_{11}$

or

```
RankOfSplitg2[{b1 }]
```

The rank is 2, and a Cartan subalgebra:

$v_1 = b_1$

$v_2 = b_4 + b_{10}$

## ■ Lie algebra so(4,4)⊕so(2,2)

Finally, we present the Lie algebra (see the references) so(4,4)⊕so(2,2). This Lie algebra is a model we have created that explains all the quarks and their corresponding antiquarks. We define a basis by making

```
θ = Table[0, {i, 1, 12}, {j, 1, 12}];
Do[m = θ; m[[p, q]] = 1;
  F_{p,q} = m, {q, 1, 12}, {p, 1, 12}];
Do[J_{p,q} = F_{p,q} - F_{q,p}; Z_{p,q} = F_{p,q} + F_{q,p},
  {p, 1, 12}, {q, 1, 12}];
bas = {Z_{1,5}, Z_{1,6}, Z_{1,7}, Z_{1,8},
    Z_{2,5}, Z_{2,6}, Z_{2,7}, Z_{2,8}, Z_{3,5},
     Z_{3,6}, Z_{3,7}, Z_{3,8}, Z_{4,5}, Z_{4,6},
    Z_{4,7}, Z_{4,8}, J_{1,2}, J_{1,3},
 J_{1,4}, J_{2,3}, J_{2,4}, J_{3,4}, J_{5,6},
    J_{5,7}, J_{5,8}, J_{6,7}, J_{6,8},
     J_{7,8}, J_{9,10}, J_{11,12}, Z_{9,11},
    Z_{9,12}, Z_{10,11}, Z_{10,12}};
n = Length[bas];
t = 12; Do[b_i = bas[[i]], {i, 1, n}];
```

$$\text{Scala[x_, y_]} := \sum_{i=1}^{t}\sum_{j=1}^{t} x[[i, j]]\, y[[i, j]];$$

$$X = \sum_{i=1}^{n} a_i\, b_i;$$

and then a generic element is

**X**

$$\begin{pmatrix}
0 & a_{17} & a_{18} & a_{19} & a_1 & a_2 & a_3 & a_4 & 0 & 0 & 0 & 0 \\
-a_{17} & 0 & a_{20} & a_{21} & a_5 & a_6 & a_7 & a_8 & 0 & 0 & 0 & 0 \\
-a_{18} & -a_{20} & 0 & a_{22} & a_9 & a_{10} & a_{11} & a_{12} & 0 & 0 & 0 & 0 \\
-a_{19} & -a_{21} & -a_{22} & 0 & a_{13} & a_{14} & a_{15} & a_{16} & 0 & 0 & 0 & 0 \\
a_1 & a_5 & a_9 & a_{13} & 0 & a_{23} & a_{24} & a_{25} & 0 & 0 & 0 & 0 \\
a_2 & a_6 & a_{10} & a_{14} & -a_{23} & 0 & a_{26} & a_{27} & 0 & 0 & 0 & 0 \\
a_3 & a_7 & a_{11} & a_{15} & -a_{24} & -a_{26} & 0 & a_{28} & 0 & 0 & 0 & 0 \\
a_4 & a_8 & a_{12} & a_{16} & -a_{25} & -a_{27} & -a_{28} & 0 & 0 & 0 & 0 & 0 \\
0 & 0 & 0 & 0 & 0 & 0 & 0 & 0 & 0 & a_{29} & a_{31} & a_{32} \\
0 & 0 & 0 & 0 & 0 & 0 & 0 & 0 & -a_{29} & 0 & a_{33} & a_{34} \\
0 & 0 & 0 & 0 & 0 & 0 & 0 & 0 & a_{31} & a_{33} & 0 & a_{30} \\
0 & 0 & 0 & 0 & 0 & 0 & 0 & 0 & a_{32} & a_{34} & -a_{30} & 0
\end{pmatrix}$$

We define then

```
RankOfModel[dat_] := Module[{},
r = Length[dat];
Do[v_k = dat[[k]], {k, 1, r}];
     If[Indepen[v, r] == 1 && Co[v, r] == 1 && Star[v, r] ==
         pset = KerAd[v_1 , b, n];
                  If[r > 1, Do[pset = Cap[pset, KerAd[v_i ,

Y = Sum[f_k Sum[pset[[k, i]] b_i, {i, 1, n}], {k, 1, Length[pset]}];

While[Maximal[v, r, b, n] == 0, r = r + 1; Addone[r, v, b,
Print["The rank is ", r, ",", " and a Cartan subalgebr
Do[comb = Solve[v_q == X];

     Print["v"_q, " = ", Sum[(a_p //. comb[[1]]) * b_p, {p, 1, n}], {q, 1,

                  Print[False]]]
```

We can now compute the rank of this Lie algebra by using two elements:

```
RankOfModel[{b1, b6}]
```

The rank is 6, and a Cartan subalgebra:

$v_1 = b_1$

$v_2 = b_6$

$v_3 = b_{34}$

$v_4 = b_{31}$

$v_5 = b_{16}$

$v_6 = b_{11}$

or with three elements

```
RankOfModel[{b1, b6, b11}]
```

The rank is 6, and a Cartan subalgebra:

$v_1 = b_1$

$v_2 = b_6$

$v_3 = b_{11}$

$v_4 = b_{34}$

$v_5 = b_{31}$

$v_6 = b_{16}$

This is the reason why this model explains six observables by acting on $\mathbb{R}^{12}$ (6 quarks + 6 antiquarks).